\documentstyle[11pt,righttag,amscd,amssymb,pstricks,verbatim]{amsart}
\newtheorem{thm}{Theorem}[section]
\newtheorem{lem}[thm]{Lemma}
\newtheorem{cor}[thm]{Corollary}
\newtheorem{prop}[thm]{Proposition}

\theoremstyle{definition}
\newtheorem{example}[thm]{Example}

\newtheorem{rem}[thm]{Remark}
\newtheorem{rems}[thm]{Remarks}

\def\lan{\langle}    \def\ran{\rangle}
\def\lr#1{\langle #1\rangle}
\def\az{\alpha}  
\def\bz{\beta}  
\def\gz{\gamma}
  \def\ooz{\Omega}
\def\sz{\sigma}

\def\vez{\varepsilon}  
\def\dz{\delta}

\def\llz{\Lambda}
\def\lz{\lambda}
\def\cm{{\cal M}}       
\def\ch{{\cal H}}     
\def\cz{{\cal Z}} \def\co{{\cal O}} \def\ce{{\cal E}}
\def\ci{{\cal I}}  
\def\ca{{\cal A}}
\def\bbn{{\mathbb N}}  \def\bbz{{\mathbb Z}}  \def\bbq{{\mathbb Q}}

 \def\vphi{\varphi}
        \def\bd{{\bf d}} 
\def\bU{{\bf U}}  
  \def\leq{\leqslant}  \def\geq{\geqslant}
\def\ge{\geqslant} \def\led{\leqslant} \def\le{\leqslant}
\def\lra{\longrightarrow}   

\def\ra{\rightarrow}
\def\hom{\mbox{\rm Hom}}

\def\ext{\mbox{\rm Ext}\,}   
\def\dim{\mbox{\rm dim}\,}
\def\udim{{\mathbf dim\,}}

  \def\End{\text{\rm End}}
\def\fkm{{\frak m}}
\def\tfu{{\tilde{\frak u}}}
\def\fkc{{\frak c}}

\begin{document}
\title[On bases of quantized enveloping algebras]%
{On bases of quantized enveloping algebras}
\author{Bangming Deng and Jie Du}
\address{Department of Mathematics, Beijing Normal University,
Beijing 100875,  China.}
\email{dengbm@@bnu.edu.cn}
\address{School of Mathematics, University of New South Wales,
Sydney 2052, Australia.} \email{j.du@@unsw.edu.au\qquad {\it Home
page}: http://www.maths.unsw.edu.au/$\sim$jied}


\thanks{Supported partially by the the NSF of China (Grant
no. 10271014), the TRAPOYP, and the Australian Research Council.}

\begin{abstract}
We give a systematic description of many monomial bases for a
given quantized enveloping algebra and of many integral monomial
bases for the associated Lusztig $\bbz[v,v^{-1}]$-form. The
relations between monomial bases, PBW bases and canonical bases
are also discussed.
\end{abstract}

\maketitle

\section{Introduction}

Let $\frak g$ be a (complex) semisimple Lie algebra and let
$\bU^+$ be the positive part of its associated quantized
enveloping algebra $\bU=\bU_v({\frak g})$ over $\bbq(v)$ with a
Drinfeld-Jimbo presentation in the generators $E_i,F_i,
K_i^{\pm1}$ $(i\in I=[1,n])$. We denote by $U^+$ the Lusztig form
of $\bU^+$, that is, $U^+$ is generated by all the divided powers
$E_i^{(m)}$ over $\cz:=\bbz[v,v^{-1}]$. Let $\ooz$ be the set of
words on the alphabet $I$ and, for $w=i_1^{e_1}i_2^{e_2}\cdots
i_m^{e_m}\in\ooz$ with $i_{j-1}\neq i_j\,\forall j$,  put
$\fkm^{(w)}=E_{i_1}^{(e_1)}\cdots E_{i_m}^{(e_m)}$. Let further
$\llz$ denote the set of all functions from the set of positive
roots of $\frak g$ to non-negative integers. The main result of
the paper is the following.

\begin{thm} Assume that $\frak g$ is simply-laced. Then there is a
partition $\ooz=\cup_{\lz\in\llz}\ooz_\lz$ such that, for any
chosen $w_\lz\in\ooz_\lz$ $(\lz\in\llz)$, the set of monomials
$\{\fkm^{(w_\lz)}\}_{\lz\in\llz}$ forms a basis for $\bU^+$.
Moreover, if all words $w_\lz$ are chosen to be distinguished,
then the set forms a $\cz$-basis for $U^+$.
\end{thm}

This work generalizes some constructions of monomial bases given
in \cite{R95} and \cite{Re1}; see \cite{DD} for a similar result
in the affine ${\frak{sl}}_n$ case. The assumption of simply-laced
types is made so that we may directly use the theory of quiver
representations, especially the theory of generic extensions
developed recently in \cite{Re}. It is natural to expect that a
similar result holds in the non-simply-laced case.

The main ingredients for the proof are Ringel's Hall algebra
theory, Reineke's monoidal structure on the set $\cm $ of
isoclasses of finite dimensional representations of a Dynkin
quiver $Q$ and the Bruhat-Chevalley type partial ordering on $\cm
$. These will be discussed separately in \S2, \S3 and \S4.
Distinguished words are introduced and investigated in \S 5 and we
shall prove the main result in \S6. As an application of the
theory, we mention an elementary construction (see \cite[\S6]{Re})
of the canonical bases for $U^+$ as the counterpart of a similar
construction for Hecke algebra in \cite{KL}. This construction
uses the same order as the one used in the geometric construction
which involves perverse sheaf and intersection cohomology
theories. Finally, more explicit results on distinguished words
are worked out for the case of type $A$ in the last section.


Throughout, $k$ denotes a finite field unless otherwise specified.
Let $q_k=|k|$. All modules are finite dimensional over $k$. If $M$
is a module, $nM$, $n\geq 0$, denotes the direct sum of n copies
of $M$. Further, by $[M]$ we denote the class of modules
isomorphic to $M$, i.e., the isoclass of $M$. For modules $M,
N_1,\cdots, N_t$, let $F^M_{N_1\cdots N_t}$ denote the number of
filtrations
$$M=M_0\supset M_1\supset \cdots \supset M_{t-1}\supset M_t=0$$
such that $M_{i-1}/M_i\cong N_i$ for all $1\le i\le t$.\\

\noindent

\section{Ringel-Hall algebras of Dynkin quivers}

Let $Q=(I,Q_1)$ be a quiver, i.e., a finite directed graph, where
$I=Q_0$ is the set of vertices $\{1,2,\cdots,n\}$ and $Q_1$ is the
set of arrows. If $\rho\in Q_1$ is an arrow from tail $i$ to head
$j$, we write $h(\rho)$ for $j$ and $t(\rho)$ for $i$. Thus we
obtain functions $h,t:Q_1\to I$. A vertex $i\in I$ is called a
sink (resp. source) if there is no arrow $\rho$ with $t(\rho)=i$
(resp. $h(\rho)=i$).

Let $kQ$ be the path algebra of $Q$. A (finite dimensional) {\it
representation} $V$ of $Q$, consisting of a set of finite
dimensional vector spaces $V_i$ for each $i\in I$ and a set of
linear transformations $V_\rho: V_{t(\rho)}\to V_{h(\rho)}$ for
each $\rho\in Q_1$, is identified with a (left) $kQ$-module. We
call $\udim V:=(\dim V_1,\cdots,\dim V_n)$ the {\it dimension
vector} of $V$ and $\ell(V):=\sum_{i=1}^n\dim V_i$ the {\it
length}\footnote{This is also called the dimension of $V$.} of
$V$. In case $Q$ contains no oriented cycles, there are exactly
$n$ pairwise non-isomorphic simple $kQ$-modules $S_1,\cdots,S_n$
corresponding bijectively to the vertices of $Q$.

{\it From now on, we assume that $Q$ is a Dynkin quiver,} that is,
a quiver whose underlying graph is a (simply laced) Dynkin graph.
By Gabriel's theorem \cite{Ga}, there is a bijection between the
set of isoclasses of indecomposable $kQ$-modules and a positive
system $\Phi^+$ of the root system $\Phi$ associated with $Q$. For
any $\bz\in\Phi^+$, let $M(\bz)=M_k(\bz)$ denote the corresponding
indecomposable $kQ$-module. By the Krull-Remak-Schmidt theorem,
every $kQ$-module $M$ is isomorphic to
$$M(\lz)=M_k(\lz):=\bigoplus_{\bz\in\Phi^+}\lz(\bz)M_k(\bz),$$
for some function $\lz: \Phi^+\to\Bbb N$. Thus, the isoclasses of
$kQ$-modules are indexed by the set
$$\Lambda=\{\lz:\Phi^+\to{\Bbb N}\}\cong{\Bbb N}^{|\Phi^+|}.$$

Further, by a result of Ringel \cite{R90}, for
$\lz,\mu_1,\cdots,\mu_m$ in $\llz$, there is a polynomial
$\vphi^\lz_{\mu_1\cdots \mu_m}(q)\in\bbz[q]$, called Hall
polynomial, such that for any finite field $k$ of $q_k$ elements
$$\vphi^\lz_{\mu_1\cdots\mu_m}(q_k)=F^{M_k(\lz)}_{M_k(\mu_1)\cdots M_k(\mu_m)}.$$

Let $\ca=\bbz[q]$ be the integral polynomial ring in the
indeterminate $q$. The generic (untwisted) Ringel-Hall algebra
$\ch=\ch_q(Q)$ of $Q$ over $\ca$ is by definition the free
$\ca$-module having basis $\{u_\lz|\lz\in\llz\}$, and satisfying
the multiplicative relations:
$$u_\mu u_\nu=\sum_{\lz\in\llz}\vphi^\lz_{\mu\nu}(q)u_\lz.$$
We sometimes write $u_\lz=u_{[M(\lz)]}$ in order to make certain
calculations in term of modules. For $i\in I$, we set
$u_i=u_{[S_i]}$. Clearly, $\ch  $ admits a natural
$\bbn^n$-grading by dimension vectors.

Following \cite{R93b}, we can twist the multiplication of the
Ringel-Hall algebra to obtain the positive part  $\bU^+$ of a
quantized enveloping algebra.

Let $\cz=\bbz[v,v^{-1}]$, where $v$ is an indeterminate with
$v^2=q$. The {\it twisted} Ringel-Hall algebra $\ch^\star=\ch_v^\star(Q)$ of
$Q$ is by definition the free $\cz$-module having basis
$\{u_\lz=u_{[M(\lz)]}|\lz\in\llz\}$ and satisfying the
multiplication rules
$$u_\mu\star u_\nu=v^{\lr{\mu,\nu}}u_\mu u_\nu=v^{\lr{\mu,\nu}}\sum_{\lz\in\llz}
\vphi^{\lz}_{\mu\nu}(v^2)u_\lz,$$ where
$\lr{\mu,\nu}=\dim_k\hom_{kQ}(M(\mu),
N(\nu))-\dim_k\ext^1_{kQ}(M(\mu), N(\nu))$ is the Euler form
associated with the quiver $Q$. Note that, if we define the
bilinear form $\lr{-,-}:\bbz^n\times\bbz^n\ra \bbz$ by
$$\lr{{\bf a}, {\bf b}}=\sum_{i\in I}a_ib_i -
\sum_{\rho\in Q_1}a_{t(\rho)} b_{h(\rho)},$$ where ${\bf a}=(a_1,
\cdots, a_n)$, ${\bf b}=(b_1, \cdots, b_n)$, then
$$\lr{\mu,\nu}=\lr{\udim M(\mu),\udim M(\nu)}.$$

For each $m\geq 1$, let
$[m]=\frac{v^m-v^{-m}}{v-v^{-1}}\;\text{and}\; [m]^!=[1][2]\cdots
[m]$. We define, for each $i\in I$, the divided powers
$$u_i^{(\star m)}=\frac{u_i^{\star m}}{[m]^!}\;\;
\text{and}\;\; E_i^{(m)}=\frac{E_i^m}{[m]^!},$$ in
$\ch^\star $ and $\bU^+$, respectively. Here $u_i^{\star
m}=\underbrace{u_i\star\cdots\star u_i}_m=v^{m\choose 2}u_i^m.$

The following result is due to Ringel \cite[\S7]{R95}.

\begin{prop} \label{lfm} The algebra $\ch^\star $ is generated by all
$u_i^{(\star m)}$ ($i\in I, m\geq1$). Moreover, there is a natural
isomorphism
$$\Psi:U^+\overset\sim\to\ch^\star ,\;E_i^{(m)}\mapsto u_i^{(\star m)},\;(i\in I,m\geq1).$$
\end{prop}
In the sequel, we shall identify $U^+$ with $\ch^\star $ under
this isomorphism.

\section{Generic extensions and the monoid $\cm $}

In this section, we collect some recent results on generic
extensions for quiver representations over an {\it algebraically
closed} field $k$.

Fix $\bd=(d_i)_i\in\bbn^n$ and define the affine space
$$R(\bd)=R(Q,\bd):=\prod_{\az\in Q_1}\hom_k(k^{d_{t(\az)}},k^{d_{h(\az)}})\cong
\prod_{\az\in Q_1} k^{d_{h(\az)}\times d_{t(\az)}}.$$ Thus, a
point $x=(x_\az)_\az$ of $R(\bd)$ determines a representation
$V(x)$ of $Q$. The algebraic group $GL(\bd)=\prod_{i=1}^n
GL_{d_i}(k)$ acts on $R(\bd)$  by conjugation
$$(g_i)_i\cdot(x_\az)_\az=(g_{h(\az)}x_\az g_{t(\az)}^{-1})_\az,$$
and the $GL(\bd)$-orbits $\co_x$ in $R(\bd)$ correspond
bijectively to the isoclasses $[V(x)]$ of representations of $Q$
with dimension vector $\bd$.

The stabilizer $GL(\bd)_x=\{g\in GL(\bd)|gx=x\}$ of $x$ is the
group of automorphisms of $M:=V(x)$ which is Zariski-open in
$\text{End}_{kQ}(M)$ and has dimension equal to
$\dim\text{End}_{kQ}(M)$. It follows that the orbit $\co_M:=\co_x$
of $M$ has dimension
$$\dim\co_M=\dim GL(\bd)
-\dim\mbox{End}_{kQ}(M).$$

\begin{lem}\label{Rei1}{\rm (\cite{Re})} Let $Q$ be a Dynkin quiver, i.e., a
disjoint union of oriented (simply-laced) Dynkin diagrams. For
$x\in R(\bd_1)$ and $y\in R(\bd_2)$, let $\ce(\co_x,\co_y)$ be the
set of all $z\in R(\bd)$ where $\bd=\bd_1+\bd_2$ such that $V(z)$
is an extension of some $M\in \co_x$ by some $N\in \co_y$. Then
$\ce(\co_x,\co_y)$ is irreducible.
\end{lem}

Given representations $M,N$ of $Q$, consider the extensions
$$0\to N\to E\to M\to 0$$
of $M$ by $N$. By the lemma, there is a unique (up to isomorphism)
such extension $G$ with $\dim\co_G$ maximal (i.e., with
$\dim\mbox{End}_{kQ}(G)$ minimal). We call $G$ the {\it generic
extension} of $M$ by $N$, denoted by $M*N$.

For two representations $M,N$, we say that $M$ degenerates to $N$,
or that $N$ is a degeneration of $M$, and write $[N]\led [M]$ (or
simply $N\led M$), if $\co_N\subseteq \overline {\co_M}$, the
closure of $\co_M$. Note that $N<M \iff \co_N\subseteq
\overline{\co_M}\backslash\co_M$.

\begin{rem}\label{REM} The relation $\le$ on the isoclasses is independent
of the field $k$. This is seen from the following equivalence
proved in \cite[Prop. 3.2]{Bo}: \begin{equation}\label{bon} N\le M
\iff \dim\hom(X,N)\ge\dim\hom(X,M),\forall X \end{equation} and
the fact that the dimension  $\dim\hom(X,Y)$ is the same over any
field. Thus, we may simply define a (characteristic-free) partial
order on $\llz$ by
$$\lz\leq\mu \iff M_k(\lz)\leq M_k(\mu).$$
for any given (algebraically closed) field $k$.
\end{rem}

The first part of the following results is well-known (see, for
example, \cite[1.1]{Bo}) and the other parts are proved in
\cite{Re}.

\begin{thm} \label{GE}
(1) If $0\to N\to E\to M\to 0$ is exact and non-split, then
$M\oplus N< E$.

(2) Let $M,N,X$ be representations of $Q$. Then $X\leq M\ast N$ if
and only if there exist $M'\leq M,N'\leq N$ such that $X$ is an
extension of $M'$ by $N'$. In particular, we have $M'\leq M,N'\leq
N\Longrightarrow M'*N'\leq M*N$.

(3) Let $\cm$ be the set of isoclasses of $kQ$-modules and define
a multiplication $*$ on $\cm$ by $[M]*[N]=[M*N]$ for any
$[M],[N]\in \cm$. Then $\cm$ is a monoid with identity $1=[0]$ and
the multiplication $*$ preserves the induced partial ordering on
$\cm$.

(4) $\cm$ is generated by the simple modules $[S_i]$, $i\in I$.
\end{thm}

Let $\Omega$ be the set of words on the alphabet
$I=\{1,\cdots,n\}$. For $w=i_1i_2\cdots i_m\in \Omega$, let
$\wp(w)\in\llz$ be the element defined by
\begin{equation}\label{wp}
[S_{i_1}]*\cdots*[S_{i_m}]=[M(\wp(w))].
\end{equation}
Thus, we obtain a map $\wp:\ooz\to\llz$. Note that by
the theorem above, $\wp$ is surjective and induces a partition of
$\ooz=\cup_{\lz\in\llz}\ooz_\lz$ with $\ooz_\lz=\wp^{-1}(\lz)$.
Each $\ooz_\lz$ is called a {\it fibre} of $\wp$.

Note from Remark \ref{REM} that, if we define $\lz*\mu$
($\lz,\mu\in\llz$) by $M(\lz*\mu)\cong M(\lz)*M(\mu)$, then the
element $\lz*\mu$ is well-defined, independent of the field $k$.
Note also that the multiplication $\ast$ on $\llz$ depends on the
orientation of $Q$.

\section{The poset $\llz$}

In this section we shall look at some properties of the poset
$(\llz,\le)$, where $\le$ is defined in Remark \ref{REM}.

For $w=i_1i_2\cdots i_m\in \Omega$ and $\lz\in\llz$, let
$\vphi_w^\lz$ denote the Hall polynomial $\varphi^\lz_{\mu_1\cdots
\mu_m}$, where $M(\mu_r)\cong S_{i_r}$. Thus, for a finite field $k$,
$$\vphi_w^\lz(q_k)=F^{M_k(\lz)}_{S_{i_1k}\cdots S_{i_mk}}$$
is the number of composition series of $M_k(\lz)$:
$$M_k(\lz)=M_0\supset M_1\supset \cdots \supset M_{m-1}\supset M_m=0$$
with $M_{j-1}/M_j\cong S_{i_jk}$. Such a composition series is
called a composition series of {\it type} $w$.

The following lemma is a
bit stronger than \cite[6.2]{DD}.

\begin{lem} \label{DOM}Let $w\in\ooz$ and $\mu\geq\lz$ in $\llz$. Then
$\vphi_w^\mu\not=0$ implies $\vphi_w^\lz\not=0$.
\end{lem}

\begin{pf} Let $w=i_1i_2\cdots i_m$ and $w'=i_2\cdots i_m$.
We apply induction on $m$. If $m=1$ then $\mu\geq\lz$ forces
$M(\mu)=M(\lz)$ and the result is clear. Assume now $m>1$. If
$\vphi_w^\mu\not=0$, then $\vphi_w^\mu(q_k)\not=0$ for some finite
field $k$. Thus, $M_k(\mu)$ has a
submodule $M_k'\cong M_k(\mu')$ which has a
composition series of type $w'$. Hence, $\vphi_{w'}^{\mu'}\not=0$
since $\vphi_{w'}^{\mu'}(q_k)\not=0$.
Base change to the algebraic closure $\bar k$ of $k$ gives an
exact sequence over $\bar k$ (We drop the subscripts $\bar k$.)
$$0\lra M'\lra M(\mu)\lra S_{i_1}\lra 0.$$
Thus,
$$M(\lz)\leq M(\mu)\leq S_{i_1}\ast M'.$$
By Theorem \ref{GE}(2), there exist modules $N', N''$ such that
$M(\lz)$ is an extension of $N'$ by $N''$ and
$N'\leq M'$, $ N''\leq S_{i_1}$. So we obtain an exact sequence (over $\bar
k$)
$$0\lra N'\overset f\lra M(\lz)\overset g\lra N''\lra 0.$$
Now the condition $N'\leq M'$ means $ \lz'\led\mu'$ where
$N'\cong M(\lz')$. Since $\vphi_{w'}^{\mu'}\neq0$, it follows from
induction that $\vphi_{w'}^{\lz'}\neq0$, that is, $N'$ has a
composition series of type $w'$.
On the other hand, since $S_{i_1}$ is simple, $ N''\leq S_{i_1}$
implies $N''\cong S_{i_1}$. Therefore, $M(\lz)$ has  a
composition series of type $w$, and consequently, $\vphi_w^\lz\not=0$.
\end{pf}

We now relate the partial order $\led$ to certain non-zero Hall
polynomials.

\begin{thm} \label{DOMI} Let $\lz,\mu\in\llz$. Then $\lz\leq \mu$
if and only if there exists a word $w\in\wp^{-1}(\mu)$ with
$\vphi_w^\lz\not=0$.
\end{thm}

\begin{pf} Suppose $\lz\leq \mu$. Since $\wp$ is surjective,
$\mu=\wp(w)$ for some $w\in\ooz$. By (\ref{wp}), we
see that $\vphi_w^{\wp(w)}\not=0$. Thus, Lemma \ref{DOM} implies
$\vphi_w^\lz\not=0$, as required.

Conversely, let $w=i_1i_2\cdots i_m\in\ooz$, $\lz\in\llz$, and suppose
$\vphi_w^\lz\not=0$. We use induction on $m$ to prove that
$\lz\leq\wp(w)$. If $m=1$, there is nothing to prove. Let
$m>1$ and $w'=i_2\cdots i_m$ and assume $\lz'\le\wp(w')$
whenever $\vphi_{w'}^{\lz'}\neq0$. Since $\vphi_w^\lz\not=0$, there is a
finite field  $k$ (of any given characteristic) such that
$\vphi_w^\lz(q_k)\not=0$. Thus, there is a submodule $M_k'$ of
$M_k(\lz)$ which has a composition series of type $w'$.
This implies $\vphi_{w'}^{\lz'}\neq0$ where $M_k(\lz')\cong M_k'$.
By induction,  we have $\lz'\leq \wp(w')$.

On the other hand, base change to the exact sequence
$$0\lra M_k'\lra M_k(\lz)\lra S_{i_1k}\lra 0$$
yields an exact sequence over $\bar k$
$$0\lra M'\lra M(\lz)\lra S_{i_1}\lra 0.$$
(Here again we dropped the subscripts $\bar k$.)
By Theorem \ref{GE}(2) we obtain that
$$M(\lz)\leq S_{i_1}\ast M(\lz')\leq S_{i_1}\ast M(\wp(w'))=M(\wp(w)).$$
Therefore, $ \lz\leq \wp(w)$.
\end{pf}

\section{Distinguished words}

Let $w=i_1i_2\cdots i_m$ be a word in $\ooz$. Then $w$ can be
uniquely expressed in the {\it tight form}
$w=j_1^{e_1}j_2^{e_2}\cdots j_t^{e_t}$, where $e_r\geq 1$, $1\leq
r\leq t$, and $j_r\not=j_{r+1}$ for $1\leq r\leq t-1$. Following
\cite[2.3]{R93}, a filtration
$$M=M_0\supset M_1\supset \cdots \supset M_{t-1}\supset M_t=0$$
of a module is called a {\it reduced} filtration of type $w$ if
$M_{r-1}/M_r\cong e_rS_{j_r}$, for all $1\leq r\leq t$. Note that
any reduced filtration of $M$ of type $w$ can be refined to a
composition series of $M$ of type $w$. Conversely, given a
composition series of $M$ of type $w$, there is a unique reduced
filtration of $M$ of type $w$ such that the given composition
series is a refinement of this reduced filtration. By
$\gz_w^\lz(q)$ we denote the Hall polynomial
$\vphi^\lz_{\mu_1\cdots\mu_t}(q)$, where $M(\mu_r)=e_rS_{j_r}$.
Thus, for any finite field $k$ of $q_k$ elements, $\gz_w^\lz(q_k)$
is the number of the reduced filtrations of $M_k(\lz)$ of type
$w$. A word $w$ is called {\it distinguished} if
$\gz_w^{\wp(w)}(q)=1$. Note that $w$ is distinguished if and only
if, for some algebraically closed field $k$, $M_k(\wp(w))$ has a
unique reduced filtration of type $w$ (cf. \cite[\S5]{DD}).

\begin{example}  Let $w={j_1}^{e_1}{j_2}^{e_2}\cdots {j_t}^{e_t}$
be in the tight form. If $j_1,\cdots,j_t$ are pairwise distinct
and satisfy
$$\ext^1_{kQ}(S_{j_r},S_{j_s})\neq 0\Longrightarrow r<s,$$
then $F^{M}_{N_1\cdots N_t}=0$ or $1$ for every $kQ$-module $M$,
where $N_r=e_rS_{j_r}$. Thus,  $w$ is distinguished.
\end{example}

Distinguished words will be used in the construction of integral
monomial bases for the Lusztig form. The following lemma shows
that these words are somehow evenly distributed.

\begin{lem} \label{df}
Each fibre of $\wp$ contains at least one distinguished word.
\end{lem}

\begin{pf} This follows directly from \cite[Lemma 4.5]{Re1}.
For completeness, we present here the construction of such
distinguished words.

By $\ci$ we denote the set of the isoclasses of indecomposable
representations of $Q$. Let $\ci_\ast$ be a {\it directed} partition of
$\ci$ (see \cite[\S4]{Re1}), that is, a partition of the set
$\ci$ into subsets $\ci_1,\cdots,\ci_m$ such that

a) $\ext^1_{kQ}(M,N)=0$ for all $M, N$ in the same part $\ci_r$,

b) $\ext^1_{kQ}(M,N)=0=\hom_{kQ}(N,M)$ if $M\in\ci_r$, $N\in\ci_s$, where
$1\leq r<s\leq m$.

Then, for each $\lz\in\llz$, we have a unique decomposition
$$M(\lz)=\bigoplus_{r=1}^m M_r,$$
where all the summands of $M_r$ belong to $\ci_r$, $1\leq r\leq
m$. Thus,
\begin{equation}\label{dfa}
\hom_{kQ}(M_r,M_s)\not=0\Longrightarrow r\leq s.
\end{equation}
Further, we order the vertices of $Q$ in a sequence
$i_1,i_2,\cdots,i_n$ such that, for
each $1<j\leq n$, $i_j$ is either a sink or an isolated vertex in
the full subquiver of $Q$ with vertices $\{i_1,\cdots,i_{j-1},
i_j\}$. Equivalently, $i_1,i_2,\cdots,i_n$ are ordered
to satisfy
\begin{equation}\label{dfb}
\ext^1_{kQ}(S_{i_j},S_{i_l})\not=0\;\Longrightarrow\;j<l.
\end{equation}
Let $\bd^{(r)}=(d^{(r)}_1,\cdots,d^{(r)}_n)=\udim M_r$, $1\leq
r\leq m$, and set
$$w_r=\underbrace{i_1\cdots i_1}_{d^{(r)}_{i_1}}\cdots\cdots
\underbrace{i_n\cdots i_n}_{d^{(r)}_{i_n}}$$ and $w_\lz=w_1\cdots
w_m\in\ooz$. Then \cite[Lemma 4.5]{Re1} implies that
$\wp(w_\lz)=\lz$ and $\gz_{w_\lz}^\lz(q)=1$, that is, $w_\lz$ is
distinguished.
\end{pf}

We call the distinguished words constructed above {\it directed}
distinguished words (with respect to the given directed partition
$\ci_\ast$).

We mention a special case of directed partitions $\ci_*$ where
each part $\ci_r$ contains only one isoclass. This case is
equivalent to ordering the indecomposable modules
$M(\bz_1),M(\bz_2),\cdots,M(\bz_\nu)$ such that
\begin{equation}\label{dfaa}
\hom_{kQ}(M(\bz_r),M(\bz_s))\not=0\Longrightarrow r\leq s.
\end{equation}
Note that monomial bases associated to these special
directed distinguished words have been constructed in \cite{L90} and \cite{R95};
see Remarks \ref{ddw} below.

The following example shows that a fibre of $\wp$ could
contain many words other than directed distinguished
ones.

\begin{example} \label{nonddw}
Let $Q$ denote the following quiver

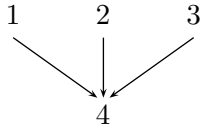
\begin{center}
\begin{pspicture}(0,0)(4,1.65)
\psset{xunit=.8cm,yunit=.8cm,linewidth=0.5pt}
\uput[u](0,1.45){1}
\uput[u](1.5,1.45){2}
\uput[u](3,1.45){3}
\uput[d](1.5,0.55){4}
\psline{->}(0,1.45)(1.4,0.45)
\psline{->}(1.5,1.45)(1.5,0.45)
\psline{->}(3,1.45)(1.6,0.45)
\end{pspicture}

\end{center}

\noindent Let $\lz\in\llz$ be such that $M(\lz)$ is the
indecomposable $kQ$-module with dimension vector $(1,1,1,2)$. Then
$\wp^{-1}(\lz)$ contains 12 words
$$\aligned
&1234^2,\, 1324^2,\, 2134^2,\, 2314^2,\, 3124^2,\, 3214^2\cr
&12434,\, 13424,\, 21434,\, 23414,\, 31424,\, 32414
\endaligned$$
which are all distinguished. From the structure of the
Auslander-Reiten quiver of $kQ$, one sees easily that the first 6
words are directed distinguished, but the last 6 words are not.
\end{example}

\section{Monomial and integral monomial bases}

For $m\geq1$, let $[\![m]\!]^!=[\![1]\!][\![2]\!]\cdots [\![m]\!]$
where ${[\![}e{]\!]}=\frac{q^e-1}{q-1}$. Then
$[\![m]\!]=v^{m-1}[m]$ and
$[\![m]\!]^!=v^{\frac{m(m-1)}{2}}[m]^!$.

\begin{lem}\label{DISL} Let $w\in\ooz$ be a word with the
tight form $j_1^{e_1}j_2^{e_2}\cdots j_t^{e_t}$. Then, for each
$\lz\in\llz$,
$$\vphi_w^\lz(q)=\gz_w^\lz(q)\prod_{r=1}^t[\![e_r]\!]^!.$$
In particular, $\vphi_w^{\wp(w)}=\prod_{r=1}^t [\![e_r]\!]^!$ if
$w$ is distinguished.
\end{lem}

\begin{pf} The result follows from the fact that the number of composition
series of $eS_i$ is $[\![e]\!]^!$ (cf. \cite[8.2]{R93b}) and the
definition of a distinguished word.
\end{pf}

To each word $ w=i_1i_2\cdots i_m\in\ooz$, we associate a monomial
$$u_w=u_{i_1}u_{i_2}\cdots u_{i_m}\in\ch.$$
Theorem \ref{DOMI} and Lemma \ref{DISL} give the following.

\begin{prop} \label{UW} For each $w\in\ooz$ with the tight
form ${j_1}^{e_1}{j_2}^{e_2}\cdots {j_t}^{e_t}$, we have
\begin{equation}\label{MtoP}
u_w=\sum_{\lz\led\wp(w)}\vphi_w^\lz(q)u_\lz =\prod_{r=1}^t
[\![e_r]\!]^! \sum_{\lz\led\wp(w)}\gz_w^\lz(q)u_\lz.
\end{equation}
Moreover, the coefficients appearing in the sum are all non-zero.
\end{prop}

We remark that this result improves \cite[Thm 1, p.96]{R95} in two
aspects. First, we generalize the formula from certain directed
distinguished words to all words; second, we
replace the lexicographical order used in \cite{R95} by the
Bruhat type partial order $\led$.

For any commutative ring $\ca'$ which is an $\ca$-algebra and any
$\ca$-module $M$, let $M_{\ca'}=\ca'\otimes_\ca M$ denote the
$\ca'$-module obtained from $M$ by base change to $\ca'$.

\begin{thm}\label{MON}  For every $\lz\in\llz$, choose an arbitrary word
$w_\lz\in\wp^{-1}(\lz)$. The set $\{u_{w_\lz}|\lz\in\llz\}$ is a
$\bbq(q)$-basis of $\ch_{\bbq(q)}$. Moreover, if all $w_\lz$
are chosen to be distinguished, then this set is an
$\ca_{(q-1)}$-basis of $\ch_{\ca_{(q-1)}}$ where
$\ca_{(q-1)}$ denotes the localization of $\ca$ at the maximal
ideal generated by $q-1$.
\end{thm}

\begin{pf} The theorem follows from Prop. 5.3 and the fact that
$\vphi_{w_\lz}^{\wp(w_\lz)}$ is invertible in $\ca_{(q-1)}$ if
$w_\lz$ is distinguished.
\end{pf}

Let ${\frak g}={\frak n}_-\oplus {\frak h}\oplus {\frak n}_+$ be
the Lie algebra over $\bbq$ of type $Q$ with generators $e_i, f_i,
h_i$. Let ${\frak U}(\frak g)$ be the universal enveloping algebra
of ${\frak g}$. We also define monomials $e_w$ similarly for
$w\in\ooz$ in ${\frak U}({\frak n}_+)$. Then, we have the
following.

\begin{cor}  For every $\lz\in\llz$, choose an arbitrary
distinguished word $w_\lz\in\wp^{-1}(\lz)$. The set
$\{e_{w_\lz}|\lz\in\llz\}$ is a $\bbq$-basis of ${\frak U}({\frak
n}_+)$.
\end{cor}

\begin{pf}
 The result
follows from  the isomorphism
 $\ch_{\ca'}/(q-1)\ch_{\ca'}\cong
{\frak U}({\frak n}_+)$,
where $\ca'=\ca_{(q-1)}$, and Theorem \ref{MON}.
\end{pf}

\noindent {\it Proof of Theorem 1.1}.  We first observe that, for
each $w=i_1i_2\cdots i_m\in\ooz$,
$$u_{i_1}\star \cdots\star u_{i_m}=v^{\vez(w)}u_w$$
where
$$\vez(w)=\sum_{1\leq r<s\leq m}\lr{\udim S_{i_r},\udim S_{i_s}}.$$
Let, for $w=j_1^{e_1}\cdots j_t^{e_t}$ in the tight form,
$$\fkm^{(w)}:=E_{j_1}^{(e_1)}\cdots E_{j_t}^{(e_t)}=(\prod_{r=1}^t [e_r]^!)^{-1}
u_{j_1}^{\star e_1}\star {\cdots} {\star} u_{j_t}^{\star e_t}.$$
Since $\prod_{r=1}^t [e_r]^!=v^{-\dz(w)}\prod_{r=1}^t [\![e_r]\!]^!$,
where
$\dz(w)=\sum_{r=1}^t \frac{e_r(e_r-1)}{2},$
it follows from Prop. \ref{UW} that
\begin{equation}\label{MtP}
\fkm^{(w)} =(\prod_{r=1}^t [\![e_r]\!]^!)^{-1}v^{\delta(w)+\vez(w)}u_w
=v^{\dz(w)+\vez(w)}\sum_{\lz\led\wp(w)}\gz_w^\lz(v^2)u_\lz.
\end{equation}
This together with Prop. \ref{lfm} and Theorem \ref{MON} implies
Theorem 1.1 with $\Omega_\lz=\wp^{-1}(\lz)$ for all $\lz\in\llz$.


\begin{rems} \label{ddw} (a) It is clear, from the definition, that
the monomial basis $\{E^{(M)}\}$ constructed in
 \cite[Thm 4.2]{Re1}
involves only directed distinguished words $w_\lz$.

(b) As a special case of \cite[Thm 4.2]{Re1}, the monomial bases
constructed in \cite[7.8]{L90} and \cite[pp.101--2]{R95}
involve only those directed distinguished words defined
with respect to the special directed partition $\ci_*$
satisfying the conditions (\ref{dfaa}) and (\ref{dfb}) in the previous section;
see \cite[Thm 1]{R95} and \cite[4.12(c),4.13]{L90}\footnote{It seems to
us that the condition (\ref{dfb})
was not explicitly stated in \cite{L90}, but was implicitly used
in \cite[7.2]{L90}.}.
\end{rems}

We now briefly look at the elementary and algebraic construction
of the canonical basis for $U^+$ (cf. \cite[\S6]{Re}). It should
be pointed out that the elementary constructions given in, e.g.,
\cite{L90}, \cite{K}, \cite{R95} and \cite{CX} used a finer order
than the one used in the geometric construction. We now use the
same order which has an algebraic interpretation (\ref{bon}).

For each $\lz\in\llz$, let
$$\tfu_\lz=v^{-\dim
M(\lz)+\dim\text{End}(M(\lz))}u_\lz.$$
Then, by Prop. \ref{lfm},
$U^+$ is $\cz$-free with basis $\ce=\{\tfu_\lz:\lz\in\llz\}$. Note
that $U^+=\oplus_{\bd}U^+_\bd$ is $\bbn I$-graded according to the
dimension vectors, and each $U^+_\bd$ is $\cz$-free with basis
$\ce\cap  U^+_\bd=\{\tfu_\lz:\lz\in\llz_\bd\}$. Clearly, each
$\llz_\bd$ together with $\led$ is a poset.

Define $\iota:U^+\to U^+$ by setting $\iota(E_i^{(m)})=E_i^{(m)}$
and $\iota(v)=v^{-1}$. Clearly, $\iota$ preserves the grading of
$U^+$. Write for any $\tfu_\lz\in U^+_\bd$
\begin{equation}\label{***}
\iota(\tfu_\mu)=\sum_{\lz}r_{\lz,\mu}\tfu_\lz.
\end{equation}
By \cite[9.10]{L90} (see \cite{Du} for more details), the
existence of the canonical bases for $U^+_\bd$ follows from the
following property on the coefficients $r_{\lz,\mu}$
\begin{equation}\label{**}
r_{\lz,\lz}=1, r_{\lz,\mu}=0\text{ unless }\lz\led \mu.
\end{equation}
 We will use (\ref{MtP}) to
derive (\ref{**}). We first calculate $\dz(w)+\vez(w)$ for
directed distinguished words (cf. \cite[Lemma, p.102]{R95}).

\begin{lem} We have for any directed distinguished word $w\in\ooz$
$$\delta(w)+\vez(w)=-\dim M(\wp(w))+\dim\text{\rm End}\,M(\wp(w)).$$
\end{lem}

\begin{pf} Let $w\in\ooz$ be a directed distinguished word. Then, by
definition, there is a directed partition $\ci_\ast$ of $\ci$ and
a $\lz\in\llz$ such that $w$ has the form $w=w_\lz=w_1\cdots w_m$
with
$$w_r=\underbrace{i_1\cdots i_1}_{d^{(r)}_{i_1}}\cdots\cdots
\underbrace{i_n\cdots i_n}_{d^{(r)}_{i_n}},$$ where
$M(\lz)=M_1\oplus M_2\oplus\cdots\oplus M_m$,
$\bd^{(r)}=(d^{(r)}_1,\cdots,d^{(r)}_n)=\udim M_r$, $1\leq r\leq
m$, and the sequence $i_1,i_2,\cdots,i_n$ of vertices are ordered
to satisfy (\ref{dfb}). Clearly, we have
$$\dz(w)=\sum_{r=1}^m\sum_{j=1}^n
\frac{d^{(r)}_{i_j}(d^{(r)}_{i_j}-1)}{2}.$$
Since $\lan\udim S_{i_j},\udim S_{i_l}\ran=0$ for $j>l$ and
$\ext^1(M_r,M_s)=0$ for all $1\leq r\le s\leq m$,
we obtain, for each $1\leq r\leq m$,
$$\aligned
\vez(w_r)&=\sum_{j=1}^n
\frac{d^{(r)}_{i_j}(d^{(r)}_{i_j}-1)}{2}\lan\udim S_{i_j},\udim
S_{i_j}\ran\cr
&+\sum_{1\leq j<l\leq n}\lan\udim
d^{(r)}_{i_j}S_{i_j},\udim d^{(r)}_{i_l}S_{i_l}\ran\cr &=\lan\udim
M_r,\udim M_r\ran -\sum_{j=1}^n\frac{(d^{(r)}_{i_j})^2}{2}
-\sum_{j=1}^n\frac{d^{(r)}_{i_j}}{2}\cr
&=\dim\End(M_r)-\sum_{j=1}^n
\frac{d^{(r)}_{i_j}(d^{(r)}_{i_j}+1)}{2}
\endaligned$$
and therefore,
$$\aligned
\vez(w)&=\sum_{r=1}^m\vez(w_r)+\sum_{1\leq r<s\leq m} \lan\udim
M_r,\udim M_s\ran\cr
&=\sum_{r=1}^m\vez(w_r)+\sum_{1\leq r<s\leq m} \dim\hom(M_r,M_s).
\endaligned$$
Noting $\hom(M_r,M_s)=0$ for $r>s$, we
finally obtain
$$\aligned
&\dz(w)+\vez(w)\cr =&\sum_{r=1}^m\dim\End(M_r)+\sum_{1\leq r<s\leq
m} \dim\hom( M_r,M_s) -\sum_{r=1}^m\sum_{j=1}^nd^{(r)}_{i_j}\cr
=&\dim\End(M(\lz))-\dim M(\lz).
\endaligned$$
This completes the proof.
\end{pf}

\begin{rem} It would be interesting to know if the lemma
holds for all distinguished words. \end{rem}

By the lemma and (\ref{MtP}), we have for a directed distinguished word $w$
\begin{equation}\label{MtooP}
\fkm^{(w)}=\tfu_{\wp(w)}+\sum_{\lz<\wp(w)}f_{\lz,\wp(w)}\tfu_\lz,
\end{equation}
where $0\neq f_{\lz,\wp(w)}\in\cz$. If we fix a representative set
$\llz'=\{w_\lz:\lz\in\llz\}$, where $w_\lz\in\ooz_\lz$, consisting of
directed distinguished words,
then the above relation implies that, for any $\lz\in\llz$,
$$\tfu_\lz\in\fkm^{(w_\lz)}+\sum_{\mu<\lz}\cz\fkm^{(w_\mu)}.$$
Restricting to $\llz_\bd$ where $\bd$ is a fixed dimension vector,
we obtain the transition matrix $(f_{\lz,\mu})_{\lz,\mu\in\llz_\bd}$. This matrix
has an inverse $(g_{\lz,\mu})_{\lz,\mu\in\llz_\bd}$ satisfying $g_{\lz,\lz}=1$
and $g_{\lz,\mu}=0$ unless $\lz\led \mu$. Thus we have
$$\tfu_\mu=\fkm^{(w_\mu)}+\sum_{\lz<\mu}g_{\lz,\mu}\fkm^{(w_\lz)}.$$
Applying $\iota$, we obtain by (\ref{MtooP})
\begin{equation}\label{MttP}
\iota(\tfu_\mu)=\fkm^{(w_\mu)}+\sum_{\lz<\mu}\bar g_{\lz,\mu}\fkm^{(w_\lz)}
=\tfu_\mu+\sum_{\lz<\mu}r_{\lz,\mu}\tfu_\lz.
\end{equation}
This proves that the coefficients in (\ref{***}) satisfy
(\ref{**}). Thus, the corresponding canonical basis
$\{\fkc_\lz\}_{\lz\in\llz}$ is uniquely defined.

\begin{rems} (a) The canonical basis defined above is the same as
Lusztig's canonical basis. This is because the basis $\ce$ is a
PBW type basis (see \cite[Thm 7]{R96}). We also note that,
as in the Hecke algebra case (\cite{KL}, \cite{KL2}), the partial order used
in this construction is the same as the one used
in the geometric construction (see \cite[\S 9]{L90}).

(b) We may also use non-directed distinguished words in the construction.
Though (\ref{MtooP}) needs to be adjusted by a power of $v$,
the relation (\ref{MttP}) will remain the same, and hence, the
canonical basis will be the same.
\end{rems}

\section{Example: the type $A$ case}

In this section, we shall give a combinatorial description of the
map $\wp:\ooz\to\llz$ for the following linear quiver:
$$Q=A_n:1\lra 2\lra\cdots \lra n-1\lra n.$$
We will also give an explicit description of the distinguished
words in this case. Since $A_n$ is a subquiver of a cyclic quiver,
results obtained below and
their proofs are similar to (or even simpler than) those given in \cite{DD}
and will be mostly omitted.

It is known that, for $1\leq i\leq j\leq n$, there is a unique (up
to isomorphism) indecomposable $kA_n$-module $M_{ij}$ with top
$S_i$ and of length $j-i+1$, and all $M_{ij}$, $1\leq i\leq j\leq
n$, form a complete set of non-isomorphic indecomposable
$kA_n$-modules. By Gabriel's theorem, each $M_{ij}$ corresponds to
a positive root $\bz_{ij}$. Thus, we have $\Phi^+=\{\bz_{ij}|1\leq
i\leq j\leq n\}$. For each map $\lz\in\llz$, we set
$\lz_{ij}=\lz(\bz_{ij})$. First, we have the following positivity
result which can be proved by counting and induction on the length
of $w$ (cf. \cite[Prop. 9.1]{DD}).

\begin{prop} For each $w\in\ooz$ and each $\lz\in\llz$,
the polynomial $\vphi_w^\lz$ lies in $\bbn[q]$.
\end{prop}

Now, for each $i\in I$, we define a map $\sz_i: \llz\ra\llz$ as
follows. For $\lz\in\llz$, if $S_{i+1}$ is not a summand of
$M(\lz)/\text{rad}M(\lz)$ (i.e., $\lz_{i+1,l}=0$, $\forall l$),
then $\sz_i\lz$ is obtained by adding 1
to $\lz_{ii}$ so that $M(\sz_i\lz)=M(\lz)\oplus S_i$; otherwise,
$\sz_i\lz$ is defined by
$$
(\sz_i\lz)_{rs}=\cases \lz_{rs}&\text{ if }\; (r,s)\neq (i,j),
(i+1,j)\cr
\lz_{ij}+1&\text{ if }\; (r,s)=(i,j)\cr
\lz_{i+1,j}-1&\text{ if }\; (r,s)=(i+1,j),\cr
\endcases
$$
where $j$ is the maximal index with $\lz_{i+1,j}\not=0$.
We have the following (cf. \cite[Prop. 3.7]{DD}).

\begin{prop}\label{SGE}  Let $i\in I$ and $\lz\in\llz$. Then
$$S_i\ast M(\lz)\cong M(\sz_i\lz).$$
Therefore, for any $w=i_1\cdots i_m\in\ooz$,
$\wp(w)=\sz_{i_1}\cdots\sz_{i_m}(0)$.
\end{prop}

Let $w=j_1^{e_1}j_2^{e_2}\cdots j_t^{e_t}\in\ooz$ be in the tight
form. For each $0\leq r\leq t$, we put
$w_r=j_{r+1}^{e_{r+1}}\cdots
j_t^{e_t}\;\;\text{and}\;\;\lz^{(r)}=\wp(w_r).$ In particular,
$w_0=w$ and $w_t=1$. Further,  for $r\geq 1$, we have
$$\lz^{(r-1)}=\wp(w_{r-1})=\underbrace{\sz_{j_r}\cdots \sz_{j_r}}_{e_r}
(\lz^{(r)}).$$

The following result gives a combinatorial description of
distinguished words (cf. \cite[5.5]{DD}).

\begin{prop}\label{DISP} Let $w=j_1^{e_1}j_2^{e_2}\cdots j_t^{e_t}\in\ooz$ and
$\lz^{(r)}$, $0\leq r\leq t$, be given as above. Then $w$ is
distinguished if and only if, for each $1\leq r\leq t$, either
$\lz^{(r)}_{j_rj}=0$ for all $j_r\le j\le n$, or
$e_r\leq \sum_{a=l_r+1}^n\lz^{(r)}_{j_r+1,a}$
where $l_r$ is the maximal index for which $\lz^{(r)}_{j_rl_r}\not=0$.
\end{prop}

\begin{pf} Using a similar argument as in \cite[Thm 5.5]{DD}, one
can show that $w$ is distinguished if and only if, for each $1\leq
r\leq t$, $M(\lz^{(r-1)})$ admits a unique submodule isomorphic to
$M(\lz^{(r)})$. However, the latter condition is equivalent to the
described combinatorial condition, as shown in \cite[Lemma
5.4]{DD}.
\end{pf}

\end{document}